\documentclass[11pt,a4]{article}
\usepackage{graphicx,amssymb,amsmath,amsthm,cite}

\DeclareMathOperator{\Span}{span}

\newtheorem{proposition}{Proposition}

\def\be{\begin{equation}}
\def\ee{\end{equation}}
\def\ben{\begin{eqnarray}}
\def\een{\end{eqnarray}}

\newcommand{\la}{\langle}
\newcommand{\ra}{\rangle}

\def\emptyy{\{0\}}
\def\h{\cal{H}}

\def\op{\hat{P}}

\newcommand{\wt}{{w}}
\newcommand{\wtk}{{w}^k}
\newcommand{\wtkp}{{w}^{k+1}}
\newcommand{\wtik}{\wt_i^k}
\newcommand{\wtikp}{\wt_i^{k+1}}

\newcommand{\fak}{{f}^k}
\newcommand{\fakp}{{f}^{k+1}}

\def\J{{\cal{J}}}

\def\jj{{\ell}}
\def\S{{\cal{S}}}
\def\SV{{\cal{V}}}
\def\SVK{{\cal{V}}_k}
\def\SVKP{{\cal{V}}_{k+1}}

\def\SW{{\cal{W}}}
\def\SWK{{\cal{W}}_k}
\def\SWC{{\cal{W}^\bot}}

\def\EVKW{\hat{E}_{\SV_k \SWC}}
\def\EVKPW{\hat{E}_{\SV_{k+1} \SWC}}

\title{Oblique Matching Pursuit}
\author{Laura Rebollo-Neira\\
Aston University, Birmingham B4 7ET, UK}
\date{\today{} at \ow}
\date{}
\begin{document}
\maketitle
\baselineskip = 2\baselineskip
\begin{abstract} 
A method for selecting a suitable subspace for
discriminating signal components through an oblique projection 
is proposed. The selection criterion is based on the consistency 
principle introduced by M. Unser and A. Aldroubi and extended 
by Y. Elder. An effective implementation 
of this principle for the purpose of subspace selection is 
achieved by updating of the dual vectors yielding the 
corresponding oblique projector. 
\end{abstract}

\section{Introduction}
Oblique projectors are of assistance to signal processing applications 
\cite{sp1,unal,sp2,sp3,sp4,yon1,sp5}, in particular 
due to their ability 
to discriminate signal components lying in different subspaces.
Thereby, as discussed in \cite{sp1}, oblique projectors are 
suitable for filtering structured noise.
Let us suppose for instance that 
a given signal $f$, represented mathematically as an element 
of a vector space $\h$, is  produced by the 
superposition of two phenomena, i.e. $f=f_1+f_2$ where 
$f_1$ belongs to a subspace $\S_1 \subset \h$ and $f_2$ 
belongs to subspace $\S_2 \subset \h$. Provided that $\S_1 \cap \S_2=
\emptyy$ we can obtain from $f$ the component $f_1$ by 
 an oblique projection onto $S_1$ along $S_2$, which 
maps $f_2$ to zero without altering $f_1$.  
The procedure is straightforward and effective if the 
corresponding subspaces $S_1$ and $S_2$, such that $S_1\cap S_2=
\emptyy$, are known \cite{sp1}. Nevertheless,
this may not be always the 
case. In this letter we address the problem of selecting 
the appropriate subspace $S_1$, from the spanning set of 
a larger subspace, in order to fulfil the condition 
$S_1\cap S_2=\emptyy$ assuming that $S_2$ is known and fixed.

Given a signal, our strategy for the selection of the 
representation subspace is in the line of 
Matching Pursuit (MP) methodologies \cite{mal1,patty,relo,tropp,swa}
and is made out of two ingredients i) the sampling/reconstruction
{\em consistency requirement} 
introduced in \cite{unal} and extended 
in \cite{yon1} ii) 
a recursive procedure for adapting  the dual vectors 
giving rise to  the corresponding oblique projector 
\cite{reo}. It will be shown 
here that the latter yields an effective implementation of a 
selection criterion that we base on the 
consistency principle. 

The letter is organized as follows: Sec~\ref{sec_2} introduces the 
general framework and discusses  
the ingredients of the approach. Namely,  the consistency 
principle and the recursive updating  of the measurement 
vectors for achieving the required oblique projection. 
The Oblique Matching Pursuit 
strategy is introduced in Sec~\ref{sec_3}. Its implementation 
is discussed in Sec~\ref{sec_4}  along with a  numerical 
example. The conclusions are drawn in Sec~\ref{sec_5}.

\section{The consistency principle and stepwise updating of 
measurement vectors} \label{sec_2} 
We represent a signal $f$ as an element of an inner product 
space that  is  assumed to be finite dimensional. 
The square norm is computed 
as $||f||^2= \la f, f\ra$, where the brackets denote the 
corresponding inner product and we define the inner product 
in such a way that if $c$ is a complex number 
$\la c f, g\ra= c^\ast \la f, g\ra$, with $c^\ast$ the complex 
conjugate of $c$. Measurements of a signal $f$ (also called samples)
will be represented as linear functionals.  
Thus a set of, say $k$, sampling vectors $w_i^{k},\,i=1,\ldots,k$ 
 provides us  with a set of $k$ measurements on
$f$ given by the inner products $\la w_i^{k},f \ra,\,i=1,\ldots,k$. 
The superscript $k$ is used to indicate that to 
reconstruct the signal we will need to modify the measurement 
vectors $w_i^{k}$ if an additional measure is considered.
From the sampling measurements we can construct an approximation 
$\fak$ of $f$ using a set of reconstruction 
vectors $v_i,\,i=1,\ldots,k$.  
The {\em consistency principle} introduced in \cite{unal} 
states that the reconstruction $\fak$ from 
$\la \wt_i^{k},f \ra,\,i=1,\ldots,k$
should be self-consistent in the sense that if the approximation is
sampled  with the same vectors the same samples should be obtained.
In other words, a consistent reconstruction must 
satisfy:
$\la \wtk_i,\fak \ra = \la \wtk_i,f \ra, \quad i=1,\ldots,k.$ 
This requirement has been  considered  further 
in \cite{yon1} where it is
proved that: {\em if the reconstruction vectors 
$v_i,\,i=1,\ldots,k$ span a subspace $\SVK$ and 
the sampling vectors 
$\wt_i^k,\,i=1,\ldots,k$ span a subspace $\SWK$  such that
its orthogonal complement $\SWC$  satisfies
$\SVK \cap \SWC = \emptyy$, then 
$\fak$ is a consistent reconstruction of $f$ if and only 
if $\fak$ is the oblique projection of $f$ onto 
$\SVK$ along $\SWC$.} We represent the corresponding
oblique projector as $\EVKW$. Hence, it is endowed with the 
following properties
i) $\EVKW^2= \EVKW$, 
ii) $\EVKW v= v, \quad \text{for any} \quad v \in \SVK$ 
iii) $\EVKW w= 0, \quad \text{for any} \quad w \in \SWC.$
Given the  conditions of the above statement, the unique consistent 
approximation of $f$ is therefore $\fak = \EVKW f$.
The oblique projector can be expressed as 
$\EVKW=\sum_{i=1}^k v_i \la  \wtk_i , \cdot \ra$
where $\la \wtk_i , \cdot \ra$ indicates that $\EVKW$ acts 
by performing 
inner products as in $\EVKW f= \sum_{i=1}^k v_i \la  \wtk_i , f \ra.$
Explicit equations for updating an oblique projector 
when a new pair of reconstruction/measurement vectors  is 
to be considered are given in \cite{reo}.
As will be discussed in the next sections, 
for the purpose of this contribution 
we can restrict the measurement vectors to be lineally 
independent. Hence the vectors $\wtkp_i$  
 yielding oblique projectors along $\SWC$ onto 
nested subspaces $\SVKP= \SVK + v_{k+1}=
\Span\{v_i\}_{i=1}^{k+1}$ can be inductively obtained as follows: 

Construct vectors $u_i= v_i - \op_{\SWC}v_i$, 
with $\op_{\SWC}$ the orthogonal projector onto 
$\SWC$. From $\wt_1^1=\frac{u_1}{||u_1||^2}$ every time 
a new vector is needed compute it, and update the 
previous ones, through the equations \cite{reo}: 
\ben
\label{eq}
\wtikp&=&\wtik - \wt_{k+1}^{k+1} \la u_{k+1}, \wtik \ra, 
\quad i=1,\ldots,k \label{wik}\\
\wt_{k+1}^{k+1}&=& \frac{q_{k+1}}{||q_{k+1}||^2},\quad 
q_{k+1}= u_{k+1}-\op_{\SWK}u_{k+1}, \label{wkkp}
\een
where $\op_{\SWK}$ is the orthogonal projector onto $\SWK= 
\Span\{u_i\}_{i=1}^k$.
It should be noticed that $\SV_{k+1}+ \SWC =\SW_{k+1} \oplus \SWC$,
with $\oplus$ indicating the orthogonal sum and $+$ the 
direct sum.

In the next section we introduce a method for stepwise 
selection of the measurement vectors aiming at finding a
subspace $\SVK$ for reconstruction  
such that $\SVK \cap \SWC= \emptyy$. 
 
\section{Oblique Matching Pursuit (OBLMP)} 
\label{sec_3}

Matching Pursuit strategies for signal representation 
 evolve by stepwise selection 
of vectors, called atoms, which are drawn from 
a large set called a dictionary. Unless the 
dictionary is orthonormal,  the seminal 
approach \cite{mal1}  does  not yield a 
stepwise reconstruction  of 
the orthogonal projection 
of the signal onto a selected subspace. A
variation of this approach, called Orthogonal Matching Pursuit (OMP)
does yield  the orthogonal projection
\cite{patty}. Such a reconstruction is therefore optimal 
in the sense of minimizing 
the norm of the approximation error. However, to render a matching 
pursuit strategy suitable for discriminating signals representing 
different phenomena,
the  approach needs to be generalized. In order to propose
the Oblique Matching Pursuit (OBLMP) method addressing this problem we 
make the following assumptions.
\begin{itemize}
\item
The subspace $\SWC$ in which the signal component to be 
filtered lies is known. 
\item
The signal we wish to filter admits a unique decomposition 
$f=f_1 + f_2$, with $f_1 \in \SVK$ and $f_2 \in \SWC$. This is 
equivalent to assuming $f \in \SVK + \SWC$ with 
$\SVK \cap \SWC= \emptyy$. 
\item
The subspace $\SVK$ can be spanned by vectors of
the dictionary in hand.
\end{itemize}
As discussed in the previous section, the reconstruction that
eliminates the signal component in $\SWC$ is 
$\fak = \EVKW f$. Our goal is to construct  the  oblique 
projector by using the appropriate dictionary vectors. We know how to 
update $\EVKW$ to $\EVKPW$ so as to account for the inclusion of an 
additional vector $v_{k+1}$. The 
question arises now as to how to select $v_{k+1}$ giving rise to 
the right subspace. We answer this 
question by recourse to the consistency principle \cite{unal,yon1}. 
Considering that at iteration  $k$ 
the approximation $\fak$ of $f$ is $\EVKW f$, let us define the 
consistency error with regard to a new  measurement $\wtkp_{k+1}$ 
as $\Delta= |\la \wtkp_{k+1}, f - \EVKW f \ra|$.
Thus to construct the approximation $\fakp= \EVKPW f$ we 
propose to select the 
measurement vector $\wtkp_{k+1}$  such that
\be
\label{cri1}
\wtkp_{k+1} = \, \text{arg\,max}_{{\jj} \in \J}\, |\la \wtkp_{{\jj}}, f - \EVKW f\ra|, 
\ee
where $\J$ is the set of indices labeling the corresponding
dictionary vectors not selected in the previous steps.
\begin{proposition}
\label{p1}
If vectors $\wtk_i,\,i=1\ldots,k$ have been selected by 
criterion \eqref{cri1} and 
$|\la\wtkp_{k+1}, f - \EVKW f\ra| \neq 0$, 
the measurement vector  
$\wtkp_{k+1}$  and the previously selected 
 vectors  $\wtk_i,\,i=1\ldots,k$ are linearly independent.
\end{proposition}
\begin{proof} Assume that, on the contrary, 
$|\la \wtkp_{k+1}, f - \EVKW f\ra| \neq 0$ 
and there exists a set of numbers $\{a_i\}_{i=1}^k$ 
 such that $\wtkp_{k+1}= \sum_{i=1}^k a_i 
\wtk_i$. Since for the previously  selected  vectors 
the consistency condition holds, i.e. 
$\la \wtk_i, f \ra = \la  \wtk_i, \EVKW f\ra,\, i=1\ldots,k$, we have 
$
|\la \wtkp_{k+1}, f - \EVKW f \ra|=  |\la
\sum_{i=1}^k a_i  \wtk_i, f - \EVKW f\ra| 
= | \sum_{i=1}^k a^\ast_i  (\la \wtk_i, f\ra - \la \wtk_i, \EVKW f\ra)|
=0.$
This contradicts our assumption, which implies that  
$\wtkp_{k+1} \ne \sum_{i=1}^k a_i \wtk_i$. 
\end{proof}
\begin{proposition} \label{p2}
All measurement vectors 
$\wtkp_{\ell}$ (c.f. eq. \eqref{cri1}) 
are orthogonal to the reconstruction 
vectors selected in previous iterations.
\end{proposition}
\begin{proof} Every $\wtkp_{\ell}$ is computed as in \eqref{wkkp} 
 and for $i=1,\ldots,k$ it is true that
$\la q_{\ell}, v_i\ra= 
\la u_{\ell}, v_i \ra - \la  \op_{\SWK} u_{\ell}, v_i\ra=
\la u_{\ell}, u_i\ra - \la u_{\ell}, \op_{\SWK} v_i\ra=
\la u_{\ell}, u_i\ra - \la u_{\ell}, u_i\ra=0$.
\end{proof}
The last proposition allows us to re-state the OBLMP 
selection criterion \eqref{cri1} as
\be
\label{cri2}
\wtkp_{k+1} = \, \text{arg\,max}_{\jj\in \J} \, |\la \wtkp_{\jj}, f \ra|.
\ee
Proposition \ref{p1} ensures that, for a given 
tolerance $\delta >0$, by stopping the 
selection process when the condition 
$\text{arg\,max}_{\jj\in \J}\, |\la \wtkp_{\jj}, f \ra| <\delta$ 
is reached, the method only selects linearly independent measurement 
vectors. Let us assume that at iteration $k+1$ the selected
indices are $\ell_1,\ldots,\ell_{k+1}$ and denote 
$u_{\ell_i}=v_{\ell_i}- \op_{\SWC} v_{\ell_i},\,i=1,\ldots,k+1$ and 
 $\wtkp_{i}, \,i=1,\ldots,k+1$ to  
the corresponding duals. Since  
$\Span\{u_{\ell_i}\}_{i=1}^{k+1}=\Span\{\wtkp_i\}_{i=1}^{k+1}$  
the fact that  $\wtkp_i,\ i=1,\ldots,{k+1}$ are 
linearly independent implies that 
 $u_{\ell_i},\,i=1,\ldots,{k+1}$ are linearly independent. 
Hence, as will be shown by the next proposition, 
at step $k+1$ the proposed selection criterion yields a 
subspace $\SVKP$ satisfying the requested 
property that $\SVKP \cap \SWC= \emptyy$. 
\begin{proposition} \label{p3}
If nonzero vectors 
$u_{\ell_i}=v_{\ell_i}- \op_{\SWC} v_{\ell_i},\,i=1,\ldots,k+1$
are linearly independent 
the only vector in $\SVKP= \Span\{v_{\ell_i}\}_{i=1}^{k+1}$ which 
is also in ${\SWC}$ is the zero vector. 
\end{proposition}
\begin{proof} Suppose that there exists $g\in \SVKP$ such that 
 $g\in \SWC$. Hence, $\op_{\SWC}g=g$ and  there exists a set 
 of numbers 
 $\{b_i\}_{i=1}^{k+1}$ to express 
 $g$ as a linear combination $g=\sum_{i=1}^{k+1} b_i v_{\ell_i}$. 
  Thus  
 $\sum_{i=1}^{k+1} b_i
 \op_{\SWC} v_{\ell_i} = \sum_{i=1}^{k+1} b_i v_{\ell_i}$, 
 which using the definition of 
 $u_{\ell_i}$ implies that $\sum_{i=1}^{k+1} b_i u_{\ell_i}=0$. 
 For nonzero linearly independent vectors this implies 
 $b_i=0, \,i=1,\ldots,k+1$ and therefore $g=0$. 
\end{proof} 
At iteration $k+1$ the selected 
indices $\ell_1,\ldots,\ell_{k+1}$ are the labels of the atoms 
 $\{v_{\ell_i}\}_{i=1}^{k+1}$ yielding the signal reconstruction 
 as given by
\be\label{rec}
f^{k+1}= \EVKPW f = \sum_{i=1}^{k+1} \la \wtkp_i, f \ra v_{\ell_i} = 
\sum_{i=1}^{k+1} c^{k+1}_i v_{\ell_i}. 
\ee
The coefficients in the last equation can be 
updated at each iteration according to 
\eqref{wik} and \eqref{wkkp}, i.e.,
\ben
c^{k+1}_{k+1}& =& \la \wtkp_{k+1}, f \ra \label{ckkp}\\
c^{k_+1}_i & =& c_i^k - c^{k+1}_{k+1} \la \wtk_i, u_{\ell_{k+1}}\ra, 
\quad
i=1,\ldots,k. \label{cik}
\een
It is appropriate to point out that these equations, 
as well as \eqref{wik} and  \eqref{wkkp}, have the identical 
form of the equations 
to modify the dual vectors and the coefficients in  
the Optimized Orthogonal Matching Pursuit Approach (OOMP) 
\cite{relo}. 
However, now the equations involve vectors of  
different nature yielding therefore a different approach. 
OOMP updating arises as the particular case, 
corresponding to $u_i \equiv v_i$, for which 
$\EVKPW  \equiv \op_{\SVKP}$. 
 Nevertheless, since the criterion for the selection process we 
have adopted here does not necessarily minimize 
the norm of the residual 
error, OOMP is not a truly particular case of the new approach. 
On the contrary, we are introducing an alternative 
selection criterion based on the consistency principle, 
which could also be considered for producing  
yet one more variation of OMP.

\section{Implementation details and numerical example} \label{sec_4}
In consistence  with the hypothesis itemized in Sec.~\ref{sec_3} 
we consider that the subspace $\SWC$ is given, i.e. 
$\{\eta_i\}_{i=1}^n$ such that $\SWC= \Span\{\eta_i\}_{i=1}^n$ 
is known.
For constructing $\op_{\SWC}$ there  are  a number of 
possibilities. 
In the example we present here the set $\{\eta_i\}_{i=1}^n$ is
 linearly dependent and we have 
used the technique for dictionary 
redundancy elimination proposed in \cite{dre}.  
 MATLAB code for its
 implementation is available at \cite{web}.
The method produces a set of orthonormal
vectors $\{\psi_i\}_{i=1}^{m}$,
 $m\le n$ that we use to construct $\op_{\SWC}= 
 \sum_{i=1}^m \psi_i \la \psi_i, \cdot \ra$. 
 
Given a dictionary  
$\{v_\ell\}_{\ell \in \J}$ 
we proceed to compute vectors 
$\{u_\ell\}_{\ell \in \J}$ as $u_\ell = v_\ell - \sum_{n=1}^m \psi_n \la \psi_n, v_\ell \ra $. 
Except for the selection criterion the next
steps parallel those for the implementation of OOMP 
but considering  the dictionary  $\{u_\ell\}_{\ell \in \J}$.
A routine for implementation of OOMP based on 
Modified Gram Smidth orthogonalization with re-orthogonalization
is also available  at \cite{web}.  With very minor changes 
 that routine can be used for the implementation of OBLMP.
The algorithm is described  below. 

Starting 
by assigning $\gamma_\ell= u_\ell, \, \ell \in \J,$ at the first
step we select the index $\ell_1$ corresponding to the index 
for which  $\la \gamma_\ell, f \ra/ || \gamma_\ell ||^2$ is maximal
and set $q_1= \gamma_{\ell_1}/ ||\gamma_{\ell_1}||, w_1^1= q_1/ ||\gamma_{\ell_1}||$ and 
$c_1^1= \la w_1^1, f \ra$. The index set $\J$ is changed to $\J= \J \setminus \ell_1$. 
At step $k+1$ the sequence $\gamma_\ell,\, \ell \in \J$ (at 
this stage $\J$ is the subset of indices not selected in the 
previous $k$ steps) is orthogonalized with respect to 
$q_{k}$ as:
$\gamma_\ell = \gamma_\ell - q_{k} \la  q_{k}, \gamma_\ell \ra $ 
and, if necessary, 
reorthogonalized with respect to  $q_1, \ldots, q_{k}$ i.e.,
$\gamma_\ell = \gamma_\ell - \sum_{j=1}^{k}  q_{j}  \la  q_{j}, \gamma_\ell \ra.$
After selecting the index $\ell_{k+1}$ as the maximizer of 
$\la \gamma_\ell, f \ra/ || \gamma_\ell ||^2$ we 
set $q_{k+1}= \gamma_{\ell_{k+1}}/ ||\gamma_{\ell_{k+1}}||, 
w_{k+1}^{k+1}= q_{k+1}/ \gamma_{\ell_{k+1}}$ and
$c_{k+1}^{k+1}= \la w_{k+1}^{k+1}, f \ra$ and 
compute $\{w_i^{k+1}\}_{i=1}^{k}$ according to 
\eqref{wik} and  $\{c_i^{k+1}\}_{i=1}^{k}$  according to \eqref{cik}.
For a given tolerance 
parameter $\delta$ the algorithm is to be stopped when  
$\la \gamma_\ell, f \ra/ || \gamma_\ell ||^2 < \delta$ for 
all $\ell \in \J$. 
The reconstructed signal is then obtained as
 in \eqref{rec}. 

We illustrate now the proposed method and its motivation
by the following example: We assume the signal space to 
be the cardinal cubic spline space with distance $0.065$
between consecutive knots, on the interval $[0,4]$. 
The background we wish to filter belongs  
to the subspace spanned by the set of functions
$\eta_i(x)={(x+1)^{-0.05i}}\,,i=1,\ldots,50,\,x\in [0,4]$. 
This set is highly redundant. A good representation
of the span can be achieved  by
just five linearly independent functions. Actually, to
avoid possible bad conditioning, we used only
three orthonormal functions for constructing $\op_{\SWC}$ and
verified a posteriori that this was enough for the backgrounds we
were dealing with. In the first test  
the dictionary is the B-spline basis on $[0,4]$. We considered  
100 signals, each of which was randomly generated as 
 linear combination  of 20 dictionary functions. One of such 
signals is plotted in the top graph of Figure~1 
added to the background.  
The functions  which  are obtained by subtracting to each  
 basis function its orthogonal projection onto $\SWC$ are not  
 exactly linearly dependent. However, the problem of constructing 
 the duals is badly conditioned. Hence, the oblique projection
 onto the whole  
 space does not yield the desired signal splitting. 
 A failed attempt to separate the signal components 
 is displayed by the broken line in the bottom  graph of Figure~1. 
 On the contrary,
 by applying the OBLMP approach, we could pick from the
 whole basis some elements spanning a subspace which 
 includes the subspace in which the signal lies. Thus,
 as depicted in the same figure, 
 the signal discrimination is successful. Equivalent results 
 were obtained for all the others signals.  In the second 
 test the dictionary spanning the identical space consists 
 of highly coherent spline atoms of  
  twice as much support as the
 corresponding basis functions
  \cite{dicacha}.  
 In this case out of 100 signals, randomly generated as
  linear combination  of 20 dictionary functions, the OBLMP approach 
  successfully split 90 of them. The  
  failures are due to the fact that, since the selection 
  process  is carried out by choosing a single atom at each step, in 
  some cases it finds the right subspace by selecting a 
  larger one which eventually includes the signal subspace.  
  The construction 
  of the duals in a larger subspace is likely 
  to become faster badly conditioned when, as in the second test, 
  the selected elements
  $u_\ell$ are 
  more coherent. On the other hand, related theoretical work  
  \cite{tropp,rel1,rel2,rel3,rel4,rel5} supports the assertion
   that a  step-wise selection approach should be expected 
  to make incorrect decisions more frequently when 
  the coherence of the dictionary is larger.

\begin{figure}[!ht]
\begin{center}
\includegraphics[width=8cm]{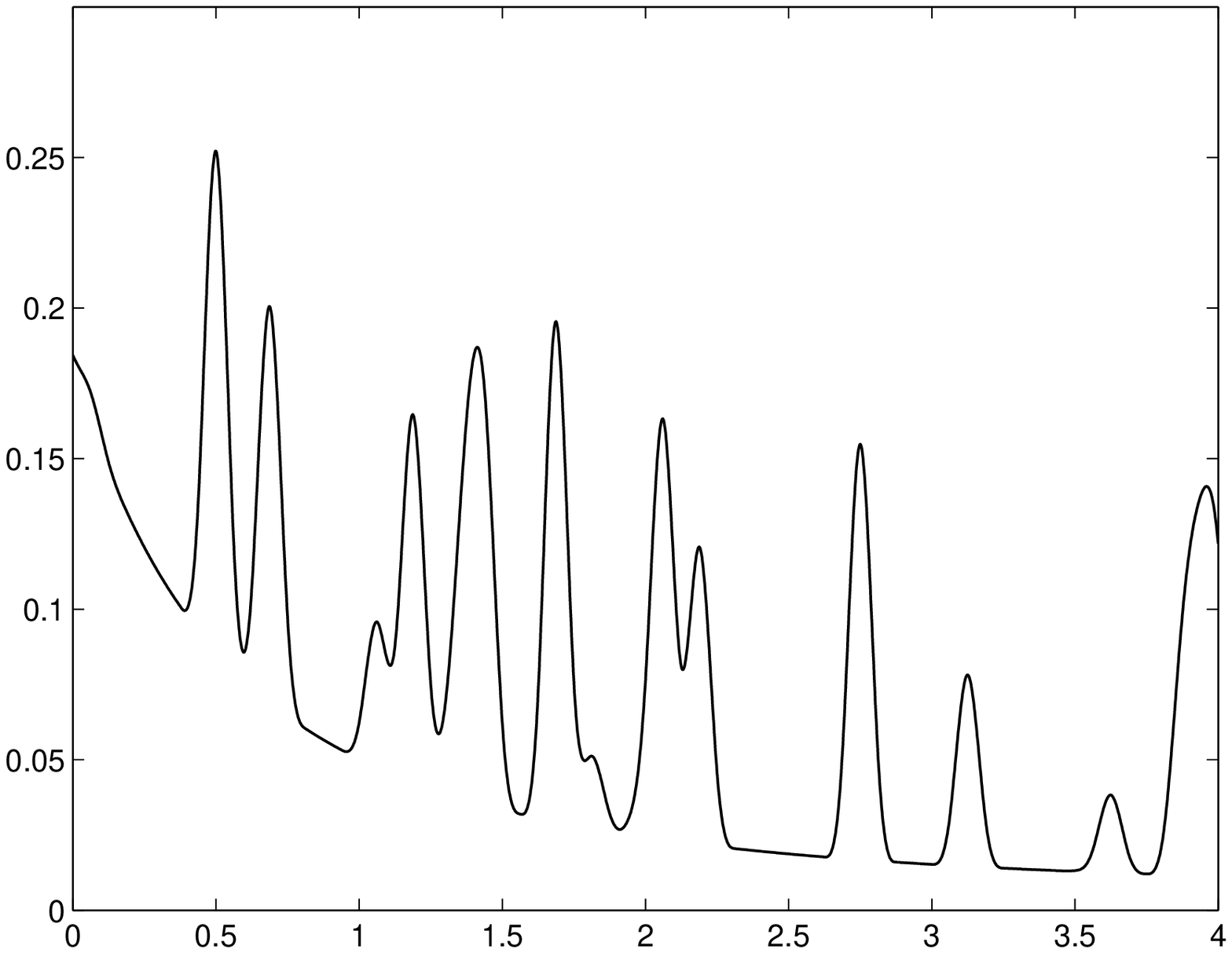}\\
\includegraphics[width=8cm]{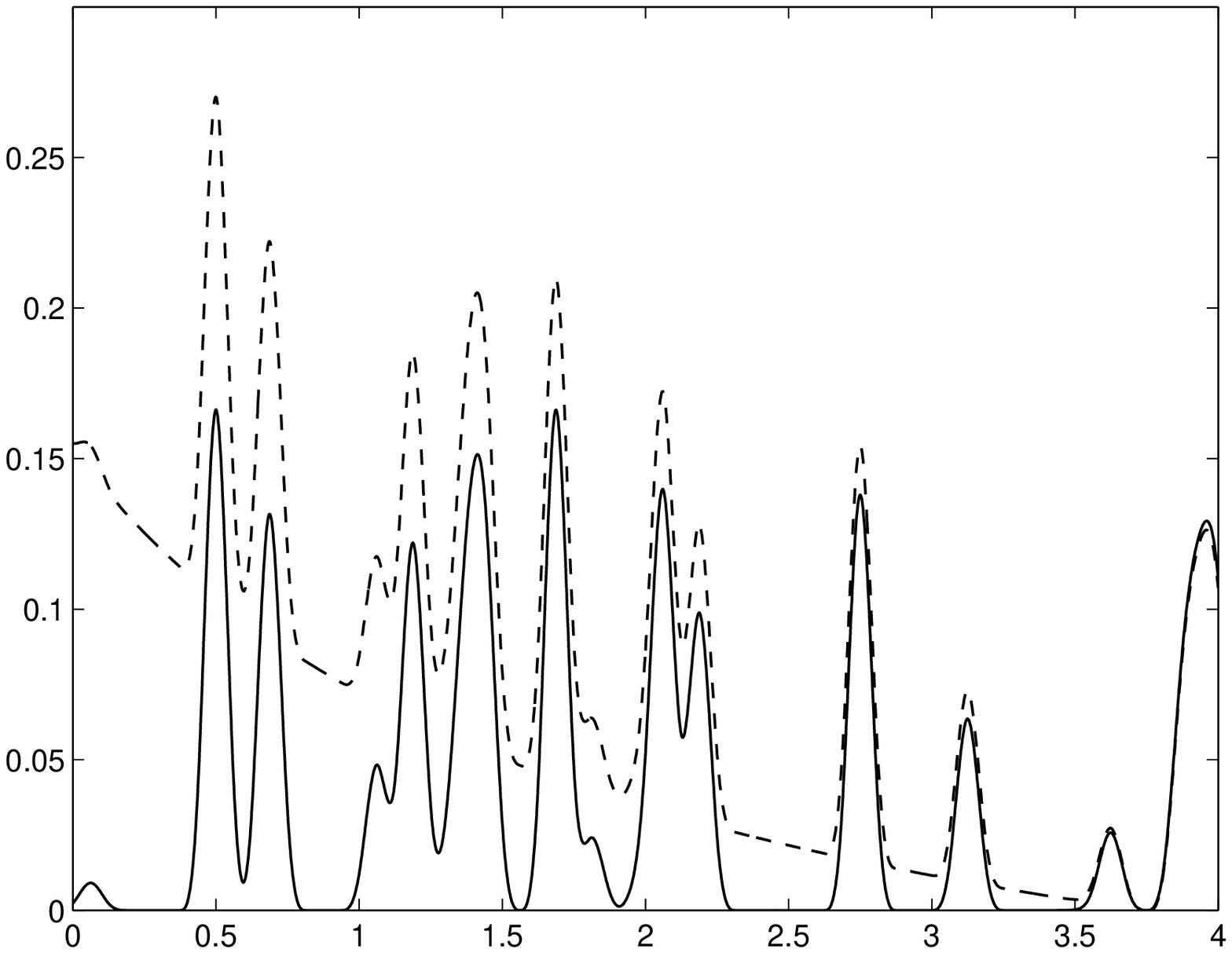}
\end{center}
\caption{The top graph shows the simulated  signal superposed on
a background belonging to the subspace $\SWC=\Span\{(x+1)^{-0.05i}\}_{i=1}^{50}$. 
The broken 
line  of the  bottom  graph depicts the result of applying the oblique 
projection onto the subspace spanned  by the whole
B-spline basis on $[0,4]$. The continuous line in the  same 
graph depicts the output of the proposed OBLMP. It reproduces the
required signal.} 
\end{figure}

\section{Conclusions} \label{sec_5}
A method, termed OBLMP, which allows for 
the selection of a suitable subspace 
for representing  one of the signal components, 
and leaving aside other components of different nature, 
has been proposed. The approach evolves by stepwise 
selection of the subspace. The selection criterion 
is based on the consistency 
requirement introduced in \cite{unal} and extended in 
\cite{yon1}.  An effective implementation is achieved by 
stepwise updating of the measurement vectors 
yielding the appropriate oblique projector \cite{reo}. 
With regard to implementation   
and complexity OBLMP is equivalent to  the
OOMP approach \cite{relo,swa}. 

Since the subspace selection is performed  by picking 
a single atom at each step, there is no guarantee that  
the required signal splitting will always be achieved. 
The success should depend on the nature of the signal components 
and the dictionaries spanning  the subspaces for 
representing them. The given examples illustrate the fact that, 
as expected, the performance of the method depends
on the coherence of the atoms resulting by subtracting 
from the dictionary atoms the orthogonal projection onto 
the background subspace. 
We hope that the results presented in this
letter will stimulate further analysis of the proposed approach. 

\subsection*{Acknowledgements}
Support from EPSRC (EP$/$D062632$/$1) is acknowledged.

\end{document}